\documentclass{amsart}
\usepackage{amsmath, amssymb, amscd, amsthm}

\newcommand{\coalg}[1]{\C{P}({#1})}

\newcommand{\lmod}[1]{{#1}\text{\rm --{\bf Mod}}}
\newcommand{\rmod}[1]{\text{{\bf Mod}--}{#1}}

\newcommand{\xra}[1]{\ensuremath{\xrightarrow{#1}}}
\newcommand{\xla}[1]{\ensuremath{\xleftarrow{#1}}}

\newcommand{\B}[1]{\ensuremath{\mathbb{#1}}}
\newcommand{\C}[1]{\ensuremath{\mathcal{#1}}}
\newcommand{\G}[1]{\ensuremath{\mathfrak{#1}}}
\newcommand{\eotimes}[1]{\ensuremath{\underset{#1}{\otimes}}}

\newcommand{\emod}[2]{\text{\bf mod}_{#1}\text{--}{#2}}
\newcommand{\eproj}[2]{\text{\bf proj}_{#1}\text{--}{#2}}
\newcommand{\efree}[2]{\text{\bf free}_{#1}\text{--}{#2}}

\newtheorem{thm}{Theorem}[section]
\newtheorem{cor}[thm]{Corollary}
\newtheorem{lem}[thm]{Lemma}
\newtheorem{prop}[thm]{Proposition}

\theoremstyle{definition}
\newtheorem{defn}[thm]{Definition}
\newtheorem{rem}[thm]{Remark}
\newtheorem{exm}[thm]{Example}

\numberwithin{equation}{section}
\title{Bivariant Hopf cyclic cohomology}
\author{Atabey Kaygun}
\email{akaygun@uwo.ca}
\author{Masoud Khalkhali}
\email{masoud@uwo.ca}
\address{Department of Mathematics, University of Western Ontario, London, Ontario, Canada}

\begin{document}

\begin{abstract}
  For module algebras and module coalgebras over an arbitrary
  bialgebra, we define two types of bivariant cyclic cohomology groups
  called bivariant Hopf cyclic cohomology and bivariant equivariant
  cyclic cohomology.  These groups are defined through an extension of
  Connes' cyclic category $\Lambda$. We show that, in the case of
  module coalgebras, bivariant Hopf cyclic cohomology specializes to
  Hopf cyclic cohomology of Connes and Moscovici and its dual version
  by fixing either one of the variables as the ground field.  We also
  prove an appropriate version of Morita invariance for both of these
  theories.
\end{abstract}

\maketitle 

\section{Introduction}

In this paper we define two types of bivariant cyclic cohomology
groups $HC^*_{\rm Hopf}(A, A'; M, M')$, and $HC^*_H(A, A'; M, M')$
where $H$ is a bialgebra, $A$ and $A'$ are $H$--module algebras and
$M$ and $M'$ are stable $H$--module/comodules.  The modules $M$ and
$M'$ appear as coefficients to twist the cyclic $H$--modules of $A$
and $A'$, respectively. We refer to the first theory as bivariant Hopf
cyclic cohomology and to the second theory as bivariant equivariant
cyclic cohomology theory. We also define bivariant groups $HC^*_{\rm
Hopf}(C, C'; M, M')$ and $HC^*_H(C, C'; M, M')$, where $C$ and $C'$
are $H$--module coalgebras.  These latter groups have much to do with
the Hopf cyclic cohomology of Hopf algebras first defined by Connes
and Moscovici \cite{ConnesMoscovici:HopfCyclicCohomology,
ConnesMoscovici:HopfCyclicCohomologyIa,
ConnesMoscovici:HopfCyclicCohomologyII} and its extensions developed
in \cite{Khalkhali:HopfCyclicHomology, Khalkhali:SaYDModules,
Khalkhali:DualCyclicHomology, Kaygun:BialgebraCyclicK}. More
precisely, we show that for $C=k$ the ground field and $M=k$ the
trivial module, we have canonical isomorphisms between $ HC^*_{\rm
Hopf}(k,C';k,M')$ and the Hopf cyclic cohomology of the $H$--module
coalgebra $C'$ with coefficients in $M'$ as defined in
\cite{Kaygun:BialgebraCyclicK}. It is shown in
\cite{Kaygun:BialgebraCyclicK} that when $H$ is a Hopf algebra and $M$
is an stable anti-Yetter--Drinfeld module this group reduces to the
Hopf cyclic cohomology of $C$ defined in
\cite{Khalkhali:HopfCyclicHomology}.  In particular for $C'=H$ a Hopf
algebra, our theory specializes to Connes-Moscovici's Hopf cyclic
cohomology for Hopf algebras.  If, on the other hand, we specialize
the second leg to $C'=k$ and $M'=k$, the groups $HC^*_{\rm Hopf}(C, k;
M, k)$ are the dual Hopf cyclic cohomology of $C$ with coefficients in
$M$. Here by duality we mean Connes' duality for the cyclic category
$\Lambda$ \cite{Connes:ExtFunctors, Connes:Book} whose relevance to
Hopf cyclic cohomology was first noticed in
\cite{Khalkhali:CyclicDuality}.

To define our bivariant groups we use the Ext bifunctor and Connes'
cyclic category $\Lambda$ in an appropriate $H$--equivariant setting.
The bivariant equivariant theory is a non-trivial extension of the
Hopf cyclic (co)homology in that there are spectral sequences
(Propositions \ref{Spectral1} and \ref{Spectral2}) converging to it
where Hopf cyclic cohomology appears in the $E_2$--term.

Our definitions immediately imply that for both theories we have
associative composition, or Yoneda products.  There is also an
analogue of Connes' periodicity $S$-operator in our context coming
from an external product: there are graded actions of the algebras
$HC^*_H(k,k;k,k)$ and $HC^*_{\rm Hopf}(k,k;k,k)$ on every bivariant
equivariant and bivariant Hopf cyclic cohomology modules
\begin{eqnarray*}
HC^p_H(k, k; k, k) \otimes HC^q_H(X, X',M, M') 
    &\to& HC^{p+q}_H(X, X'; M, M')\\
HC^p_{\rm Hopf} (k, k; k, k) \otimes HC^q_{\rm Hopf} (X, X', M, M') 
    &\to& HC^{p+q}_{\rm Hopf} (X, X', M, M')
\end{eqnarray*}
where the algebra $HC^*_{\rm Hopf}(k, k; k, k)$ is the $k$--algebra
generated by the $S$--operator.  On the other hand, the algebra
$HC^*_H(k,k;k,k)$ is a combination of ${\rm Ext}^*_H(k,k)$ and the
$k$--algebra generated by the $S$--operator via a spectral sequence of
graded $k$--algebras (Corollary~\ref{BottAlgebra}).  Thus, one can see
that bivariant equivariant cyclic cohomology is a natural extension of
both cyclic (co)homology and ordinary (co)homology of $H$--modules in
a unified theory. 

The original motivation to develop a bivariant
cyclic cohomology \cite{
JonesKassel:BivariantCyclicTheory, Nistor:BivariantChernConnes, Cuntz:BivariantKandCyclicTheories},
first defined by Connes in \cite{Connes:ExtFunctors,
Connes:NonCommutativeGeometry} as an Ext group, was to define a receptacle for a
bivariant Chern-Connes character defined on (smooth) cycles of
KK-theory. There is a similar question in our Hopf setting as well but
we will not address it in this paper.

Here is a plan of this paper.  In Section~\ref{Definitions} we define
the category of (co)cyclic $H$--modules and give a reinterpretation of
these modules, and their variations, as left and right modules over a short list
of algebras all of which can be defined in terms of a single large
algebra $\coalg{H}$.  In Section~\ref{HomologicalAlgebra} we prove
some results in homological algebra for the algebra $\coalg{H}$ and
its modules we defined in Section~\ref{Definitions}.  In
Section~\ref{ModuleAlgebra} we define a cyclic $H$--module associated
with an $H$--module algebra with coefficients in a stable
$H$--module/comodule.  Using this object we define bivariant Hopf and
bivariant equivariant cyclic homology of a pair of $H$--module
algebras with coefficients in an arbitrary pair of stable
$H$--module/comodules. By using the results of
Section~\ref{HomologicalAlgebra}, we also show how our theories relate
to the ordinary cyclic homology and cohomology of algebras, in
combination with the cohomology of these algebras viewed simply as
$H$--modules.  In Section~\ref{ModuleCoalgebra} we develop the same
theory for module coalgebras and we also investigate the connections
of these bivariant theories with Hopf cyclic cohomology.  In
Section~\ref{MoritaInvariance} we develop the notion of
$H$--categories and a cyclic homology theory for $H$--categories.  We
show that (co)cyclic $H$--module associated with a module (co)algebra
can also be interpreted as an $H$--categorical invariant, which leads
us to the Morita invariance.

Throughout this paper we assume $k$ is a field and $H$ is an
associative/coassociative, unital/counital bialgebra, or a Hopf
algebra with an invertible antipode whenever it is necessary. By a
(say, left-left) $H$--module/comodule we mean a left $H$--module which
is also a left $H$--comodule with no compatibility assumption between
action and coaction.

\section{The category of cyclic $H$--modules and its variants}
\label{Definitions}

In this section we define the categories of cyclic and cocyclic
$H$--modules where $H$ is a bialgebra. We also define closely related
categories of para-(co)cyclic and pseudo-para-(co)cyclic
$H$--modules. We show that when $H$ is a Hopf algebra, Connes'
fundamental isomorphism between the cyclic category and its dual
\cite{Connes:ExtFunctors}, can be extended to all of the above
categories. This plays an important role in our definition of
bivariant Hopf cyclic groups.

We denote Connes's cyclic category by $\Lambda$. Recall from
\cite{Connes:ExtFunctors, Connes:Book, Loday:CyclicHomology} that a
cyclic (resp. cocyclic) object in a category $\mathcal{C}$ is a
contravariant (resp. covariant) functor $X_\bullet:\Lambda \to
\mathcal{C}$. A cyclic (resp. cocyclic) $k$--module is simply a cyclic
(resp. cocyclic) object in the category of $k$--modules. Given a
cyclic module $X_{\bullet}$ we denote its cyclic homology groups by
$HC_{*}(X_{\bullet})$. Similarly we write $HC^{*}(X_{\bullet})$ to
denote the cyclic cohomology of a cocyclic module $X_{\bullet}$.

A remarkable property of Connes' cyclic category is its self
duality. It is shown in \cite{Connes:ExtFunctors} that there is a
natural isomorphism of categories $\Lambda\xra{}\Lambda^{op}$. This
fact plays almost no role in the cyclic homology of algebras or
coalgebras, but it is of considerable importance in Hopf cyclic
cohomology as it was first observed in \cite{Khalkhali:CyclicDuality}.
In the sequel the cyclic duality unavoidably manifests itself in the
Specialization Theorem (Theorem~\ref{Specialization}).

\begin{defn} 
  Let $H$ be a bialgebra. A (co)cyclic module $X_\bullet$ is called a
  (co)cyclic $H$--module if (i) for each $n\geq 0$, the module $X_n$
  is a (left) right $H$--module and (ii) all the (co)cyclic structure
  morphisms commute with the action of $H$.  A para-(co)cyclic
  $H$--module is the same as a (co)cyclic $H$--module except that the
  cyclic operators $\tau_n$ do not necessarily satisfy
  $\tau_n^{n+1}=id_n$ for $n\geq 0$.  We may even drop the condition
  that $\tau_n$ is invertible, if necessary.
\end{defn}

To define our bivariant groups we have to reinterpret the (co)cyclic
$H$--modules and its relatives as modules over certain (non-unital and
non-counital) bialgebras defined below.  Here we define them as
algebras but we are going to prove in Proposition~\ref{BigBialgebra}
that they really are non-unital and non-counital bialgebras.
\begin{defn}
  Let $H[\Lambda_\B{N}]$ be the algebra generated by $H$--linear
  combinations of the symbols $\partial^n_i$, $\sigma^m_j$ and
  $\tau_n^\ell$ for $0\leq n$, $0\leq m$, $0\leq i\leq n+1$, $0\leq
  j\leq m$ and $\ell\in\B{N}$ satisfying
  \begin{align}
    \partial^{n+1}_i\partial^n_j = & \partial^{n+1}_{j+1}\partial^n_i \text{ \ \ and \ \ }
    \sigma^{n-1}_j\sigma^n_i = \sigma^{n-1}_i\sigma^n_{j+1} \text{ for } i\leq j \ \  \text{ and }
    \tau_n^s\tau_n^t = \tau_n^{s+t} \text{ for } s,t\in\B{N} \label{Category1}\\
    \sigma^n_i\partial^n_i = & \sigma^n_i\partial^n_{i+1} = \tau_n^0 \text{ \ \ and\ \ }
    \partial^n_i\sigma^n_j 
      = \begin{cases}
	\sigma^{n+1}_{j+1}\partial^{n+1}_i & \text{ if }i\leq j\\
	\sigma^{n+1}_j\partial^{n+1}_{i+1} & \text{ if }i>j
	\end{cases} \label{Category2}\\
    \partial^n_i\tau_n^j = & \tau_{n+1}^{i+p}\partial^n_q 
        \text{ where $(i+j)=(n+1)p+q$ with $0\leq q\leq n$} \label{Category3}\\
    \tau_n^i\sigma^n_j = & \sigma^n_q\tau_{n+1}^{i+p} 
        \text{ where $(-i+j)=(n+1)(-p)+q$ with $0\leq q\leq n$}\label{Category4}
  \end{align}
  All other products between the generators are $0$ and the product is
  extended $H$--linearly.  We also define another algebra
  $H[\Lambda_\B{Z}]$ where we allow $\tau_n^\ell$ with
  $\ell\in\B{Z}$. Define $H[\C{D}]$ as the subalgebra generated by
  $\partial^m_j$ and $H[\C{D}_\B{N}]$ as the subalgebra generated by
  $\tau_n^i$ and $\partial^m_j$ for all possible $n,m,j$ and $i\geq
  0$.  We finally define $H[\Lambda_+]$ as the subalgebra generated by
  $\sigma^n_j$ and $\partial^n_j$ where we only require $0\leq j\leq
  n$ for any $n\geq 0$.
\end{defn}

\begin{defn}
  Define a $k$--algebra $\coalg{H}$ as an amalgamated product of
  $k$--algebras
  \begin{align*}
    \coalg{H}:= k[\Lambda_\B{Z}]\underset{k[\Lambda_+]}{*}H[\Lambda_+]
  \end{align*}
  It is the $k$--algebra generated by the same generators and
  relations as $H[\Lambda_\B{N}]$ however we drop the relations
  $[x,\tau_n^i]=0$ for any $x\in H$ and for all possible $n\geq 0$ and
  $i$.  Any left $\coalg{H}$--module is called a pseudo-para-cocyclic
  $H$--module and any right $\coalg{H}$--module is called a
  pseudo-para-cyclic $H$--module.  We define another algebra
  $\coalg{H}'$ by $k[\Lambda_\B{N}]\underset{k[\Lambda_+]}{*}
  H[\Lambda_+]$ where we drop the condition that $\tau_n$ is
  invertible for $n\geq 0$.
\end{defn}

\begin{rem}
  The difference between $H[\Lambda_\B{Z}]$ and $\coalg{H}$ is that
  for para-(co)cyclic modules (i.e. modules of $H[\Lambda_\B{Z}]$),
  the action of elements of $H$ and the action of the generators
  $\tau_n^\ell$ do necessarily commute while for
  pseudo-para-(co)cyclic (i.e. modules of $\coalg{H}$) this is not
  required.
\end{rem}

\begin{defn}
  Let $H[\Lambda]$ be the quotient of $H[\Lambda_\B{N}]$ by the right
  (which is also a bilateral) ideal generated by elements of the form
  $(\tau_n^{n+1}-\tau_n^0)$ for all $x\in H$ and for all possible $n$
  and $i$.  Define also $H[\C{D}_C]$ as the image of $H[\C{D}_\B{N}]$
  under this quotient.
\end{defn}

\begin{defn}
  Let $\C{A}$ be one of $\coalg{H}$, $\coalg{H}'$, $H[\Lambda_\B{N}]$,
  $H[\Lambda_\B{Z}]$, $H[\Lambda]$, $H[\C{D}_\B{N}]$ and $H[\C{D}_C]$.
  A left (or right) $\C{A}$--module $X_\bullet$ is called faithful if
  for every $x\in X_\bullet$ there exists a finite number of elements $y_i\in
  X_\bullet$ and integers $n_i\in\B{N}$ such that $x=\sum_i \tau_{n_i}^0 y_i$
  (resp. $x = \sum_i y_i\tau_{n_i}^0$).
\end{defn}

Throughout the paper, we will assume that all of our modules over the 
algebras $\C{A}$ considered  above are faithful.

\begin{lem}
  Let $\C{A}$ be one of the algebras $H[\Lambda_\B{Z}]$, $H[\Lambda]$,
  $H[\C{D}_\B{N}]$, $H[\C{D}_C]$.  Then the category of faithful right
  (left) $\C{A}$--modules is isomorphic to the category of
  para-(co)cyclic, (co)cyclic, pre-para-(co)cyclic and pre-(co)cyclic
  $H$--modules, respectively.
\end{lem}

\begin{proof}
  We will give the proof for pre-para-cocyclic case but the proof for
  the other cases are very similar.  Assume $X_\bullet$ is a
  pre-para-cocyclic $H$--module.  Then $X_\bullet$ is a
  $\B{N}$--graded $H$--module with structure morphisms $\partial_j^n$,
  and $\tau_n^\ell$ where $\ell,n\in\B{N}$ and $0\leq j\leq n+1$ which
  satisfy the conditions stated in Equations~(\ref{Category1}) through
  (\ref{Category4}) and the condition that $\tau_n^0(x)=x$ for any
  $x\in X_n$.  In other words $X_\bullet$ is a faithful left
  $H[\C{D}_\B{N}]$--module.  Conversely, assume $X_\bullet$ is a
  faithful left $H[\C{D}_\B{N}]$--module.  Define $X_n$ as the
  submodule of $X_\bullet$ consisting of elements $x$ such that
  $x=\tau_n^0(y)$ for some $y\in X_\bullet$.  Note that since
  $\tau_n^0$ is idempotent $\tau_n^0(x)=x$ for any $x\in X_n$.
  Moreover, since $X_\bullet$ is faithful $X_\bullet=\bigoplus_n X_n$
  by definition.  The actions of the generators $\partial^n_j$ define
  a pre-cosimplicial structure on $X_\bullet$.  For every $n\geq 0$ we
  also have an action of $\B{N}$ on $X_n$ via $\tau_n^\ell$.  We must
  check that combination of these structure morphisms do really define
  a pre-para-cocyclic structure.  Consider Equation~(\ref{Category3})
  for $i=n+1$ and $j=0$ and we see that since $i+j=n+1=(n+1)1+0$ we
  must have
  \begin{align*}
    \partial^n_{n+1} = \tau_{n+1}\partial^n_0
  \end{align*}
  Moreover, we observe that since $j+1=(n+1)0+j+1$ for $0\leq j\leq
  n-1$ we get
  \begin{align*}
    \partial^n_j\tau_n 
    = \begin{cases}
      \tau_{n+1}\partial^n_{i+1} & \text{ if } 0\leq j\leq n-1\\
      \tau_{n+1}^2\partial^n_0   & \text{ if } j=n
      \end{cases}
    = \tau_{n+1}\partial^n_{j+1}
  \end{align*}
  i.e. $X_\bullet$ is a pre-para-cyclic $H$--module.
\end{proof}

From now on we will use the terms left (right) $\coalg{H}$--module,
$H[\Lambda_\B{Z}]$--module, $H[\Lambda]$--module,
$H[\C{D}_\B{N}]$--module and $H[\C{D}_C]$--modules and
pseudo-para-(co)cyclic, para-(co)cyclic, (co)cyclic,
pre-para-(co)cyclic and pre-(co)cyclic $H$--module interchangeably.

The algebra $\coalg{H}$ plays an important role for us since all the
algebras we described above can be obtained from $\coalg{H}$.  Almost
all important properties of $\coalg{H}$ will descend on the rest of
the algebras.

\begin{prop}
  Let $\C{A}$ be one of $\coalg{H}$, $H[\Lambda_\B{Z}]$, $H[\Lambda]$
  and $H[\C{D}_C]$.  If $H$ is a Hopf algebra with an invertible
  antipode then $\C{A}$ is isomorphic to its opposite algebra via an
  isomorphism $\C{A}\xra{\gamma}\C{A}^{op}$.
\end{prop}

\begin{proof}
  We will give the isomorphism on $\coalg{H}$.  The isomorphism for
  the other cases is obtained from this isomorphism.  We define
  $\coalg{H}\xra{\gamma}\coalg{H}^{op}$ by defining it on generators
  as
  \begin{align*}
    \gamma(\partial^n_j) = & \sigma^n_j, \ \ \gamma(\sigma^n_j) = \partial^n_j, \ \ 
    \gamma(\tau_n)=\tau_n^{-1} \text{ \ \ and\ \ } \gamma(h) = S(h)
  \end{align*}
  for $h\in H$, $0\leq j\leq n+1$.  Note that since
  $\partial^n_{n+1}=\tau_{n+1}\partial^n_0$, one does not have a
  problem of defining $\sigma^n_{n+1} := \gamma(\partial^n_{n+1}):= 
  \sigma^n_0\tau_{n+1}^{-1}$.  It is routine to check that $\gamma$ is
  an isomorphisms of algebras, since $S$ is invertible.
\end{proof}

\begin{thm}\label{CyclicDuality}
  Let $H$ be a Hopf algebra with an invertible antipode.  Then the
  categories of (pseudo-)(pre-)(para-)cocyclic and
  (pseudo-)(pre-)(para-)cyclic $H$--modules are isomorphic.
\end{thm}

\begin{proof}
  We are going to give the proof for pseudo-para-cyclic modules.  The
  proof for the other cases is similar.  We define a functor
  \begin{align*}
    (\ \cdot\ )^\vee:\lmod{\coalg{H}}\xra{}\rmod{\coalg{H}}
  \end{align*}
  as follows: for every pseudo-para-cocyclic $H$--module $Y_\bullet$,
  the module $Y^\vee_\bullet$ is the same as $Y_\bullet$.  However,
  for every $\Psi\in\coalg{H}$ and $y\in Y_\bullet$, we define
  $y\Psi:=\gamma(\Psi)y$.  In the opposite direction,
  $(\ \cdot\ )^\vee:\rmod{\coalg{H}}\xra{}\lmod{\coalg{H}}$ is defined
  similarly.  However, given any para-cyclic $H$--module $X_\bullet$
  and $x\in X^\vee_\bullet$ and $\Psi\in H[\Lambda_\B{Z}]$ we define
  $\Psi x := x\gamma^{-1}(\Psi)$.  This way, one can see that
  $(Z^\vee_\bullet)^\vee$ is $Z_\bullet$ itself for any
  pseudo-para-(co)cyclic module $Z_\bullet$.
\end{proof}

\begin{cor}\label{IsomorphicHom}
  Let $\C{A}$ be one of $\coalg{H}$, $H[\Lambda_\B{Z}]$, $H[\Lambda]$,
  $H[\C{D}_C]$.  Let $X_\bullet$ and $Y_\bullet$ be two right
  $\C{A}$--modules.  Then one has an isomorphism of
  $k$--modules
  \begin{align*}
    {\rm Hom}_\C{A}(X_\bullet,Y_\bullet)
    \cong {\rm Hom}_\C{A}(X^\vee_\bullet,Y^\vee_\bullet)
  \end{align*}
\end{cor}

The functor we defined in Theorem~\ref{CyclicDuality} extends the
cyclic duality functor defined by Connes \cite{Connes:ExtFunctors}.
As is observed in \cite{Khalkhali:CyclicDuality} all the classical
examples coming from the cyclic and cocyclic modules of algebras and
coalgebras are all `one-sided'.  That is, if $X$ is a (co)algebra and
if $C_*(X)$ is the ordinary (co)cyclic module associated with $X$,
then the cyclic dual $C_*(X)^\vee$ has trivial homology.  However, for
an arbitrary (co)cyclic $H$--module this need not be the case.  For
example, for a Lie algebra $\G{g}$, both $\C{C}_\bullet(U(\G{g}),k)$
(Definition~\ref{HopfCyclicModule}) and
$\C{C}_\bullet(U(\G{g}),k)^\vee$ are homologically non-trivial. For
example their periodic cyclic (co)homologies are computed in
\cite{ConnesMoscovici:HopfCyclicCohomology} and
\cite{Khalkhali:DualCyclicHomology} respectively, and shown to be both
isomorphic to the Lie algebra homology of $\G{g}$ with trivial
coefficients.

\section{(Co)Homology of (co)cyclic $H$--modules}
\label{HomologicalAlgebra}

Our goal in this section is to extend some results of
\cite{Connes:ExtFunctors} from the category of (co)cyclic $k$--modules
to the category of (co)cyclic $H$--modules.

Let $\C{A}$ be one of the algebras $\coalg{H}$, $H[\Lambda_\B{Z}]$,
$H[\Lambda]$, $H[\C{D}_C]$.  We consider the $\B{N}$--graded
$k$--module $k_\bullet = \bigoplus_{n\geq 0} k $ as a left
$\C{A}$--module by letting $H$ act by the counit $H\xra{\varepsilon}k$
and by identity for the rest of the generators.

\begin{prop}\label{CyclicHomology}
  Let $H=k$.  Then for any cyclic module $X_\bullet$ we have natural
  isomorphisms
  \begin{align*}
    {\rm Ext}_{k[\Lambda]}^*(X_\bullet, k^\vee_\bullet)\cong HC^*(X_\bullet) & & 
    {\rm Tor}^{k[\Lambda]}_*(X_\bullet,k_\bullet)\cong HC_*(X_\bullet) 
  \end{align*}
\end{prop}

\begin{proof}
  In \cite[Section 4]{Connes:ExtFunctors} a projective biresolution of
  $k_\bullet^\vee$ ($k^\sharp$ in Connes' notation) is developed.
  Connes uses this particular resolution to obtain the first part of
  the Proposition.  The second part follows immediately.
\end{proof}

\begin{lem}\label{Mirror}
  Assume we have two left $k[\Lambda]$--modules, i.e. cocyclic
  $k$--modules, $X_\bullet$ and $Y_\bullet$.  Then one has a natural
  isomorphism of $k$--modules of the form
  \begin{align*}
    X_\bullet^\vee\eotimes{k[\Lambda]}Y_\bullet 
    \xra{t} Y_\bullet^\vee\eotimes{k[\Lambda]}X_\bullet
  \end{align*}
\end{lem}

\begin{proof}
  The isomorphism is defined by the transposition $t(x\otimes y) =
  (y\otimes x)$ for any $x\otimes y$ in
  $X_\bullet^\vee\eotimes{k[\Lambda]}Y_\bullet$.  However, we must
  prove that the morphism $t$ is well-defined.  For that purpose, we
  recall that $X_\bullet^\vee\eotimes{k[\Lambda]}Y_\bullet$ is a
  quotient of $X_\bullet^\vee\otimes Y_\bullet$ by the $k$--submodule
  generated by elements of the form $(x\Psi\otimes y) - (x\otimes\Psi
  y)$ where $x\in X_\bullet^\vee$, $y\in Y_\bullet$ and $\Psi\in
  k[\Lambda]$.  However, recall from Theorem~\ref{CyclicDuality} that
  $x\Psi:=\gamma(\Psi)x$ therefore
  \begin{align*}
    t((x\Psi\otimes y) - (x\otimes\Psi y))
     = (y\otimes\gamma(\Psi)x) - (y\gamma^{-1}(\Psi)\otimes x)
  \end{align*}
  However, $H=k$, $k[\Lambda]\xra{\gamma}k[\Lambda]^{op}$ is an
  involution.  The result follows.
\end{proof}

\begin{rem}
  In the general situation when $H\neq k$, Lemma~\ref{Mirror} would
  work only if one has (i) $S^2=id$ on $H$ or (ii) one of the modules
  involved has the property that $\gamma^2(\Psi)z=\Psi z$ for any
  $\Psi$.  $k_\bullet$ and $k_\bullet^\vee$ are two important examples
  of such $H[\Lambda]$--modules while $H_\bullet$ and $H_\bullet^\vee$
  are not unless $S^2=id$.
\end{rem}

\begin{prop}\label{CocyclicHomology}
  Let $H=k$.  Then for any cocyclic module $Y_\bullet$ we have natural
  isomorphisms
  \begin{align*}
    {\rm Tor}^{k[\Lambda]}_*(k^\vee_\bullet,Y_\bullet)\cong HC_*(Y^\vee_\bullet) & & 
    {\rm Ext}_{k[\Lambda]}^*(Y_\bullet,k_\bullet) \cong HC^*(Y^\vee_\bullet)
  \end{align*}
    Moreover,
  \begin{align*}
    {\rm Ext}_{k[\Lambda]}^*(k_\bullet,Y_\bullet)\cong HC^*(Y_\bullet)
  \end{align*}
\end{prop}

\begin{proof}
  The first assertion easily follows from Lemma~\ref{Mirror}.  For the
  second assertion, we observe that $(\ \cdot\ )^\vee$ is an isomorphism
  of categories which means one has a natural isomorphism of functors
  (Corollary~\ref{IsomorphicHom})
  \begin{align*}
    {\rm Hom}_{k[\Lambda]}(\ \cdot\ , k_\bullet)\cong
    {\rm Hom}_{k[\Lambda]}((\ \cdot\ )^\vee, k_\bullet^\vee)\cong
    {\rm Hom}_k\big((\ \cdot\ )^\vee\eotimes{k[\Lambda]}k_\bullet, k\big)
  \end{align*}
  by using the fact that for the trivial cocyclic module $k_\bullet$,
  the $k$--vector space dual and the cyclic dual are isomorphic.  Then
  we use Lemma~\ref{Mirror} again.  The last part of our assertion
  uses the projective resolution $C^{*,*}$ of $k_\bullet^\vee$
  developed in \cite[Section 4]{Connes:ExtFunctors} which actually is
  a double complex.  Note that since $(\ \cdot\ )^\vee$ provides an
  isomorphism of categories, $(C^{*,*})^\vee$ is a projective
  resolution of $k_\bullet$.  We observe that for this specific
  resolution, the homology of the double complex ${\rm
  Hom}_{k[\Lambda]}((C^{*,*})^\vee,Y_\bullet)$ is exactly the cocyclic
  double complex of the cocyclic $k$--module $Y_\bullet$.
\end{proof}

\begin{prop}
  Let $X_\bullet$ be a cyclic $H$--module.  Then there are three
  spectral sequences with 
  \begin{align*}
    E^2_{p,q}=HC^p({\rm Tor}^H_q(X_\bullet,k)) & & 
    'E^2_{p,q}= {\rm Hom}_k({\rm Tor}^H_p(HC^q(X_\bullet),k),k) & & 
    ''E^2_{p,q} = {\rm Ext}_H^p(HC_q(X_\bullet),k)
  \end{align*}
  converging to ${\rm Ext}_{H[\Lambda]}^*(X_\bullet,k_\bullet^\vee)$.
\end{prop}

\begin{proof}
  Note that one can write the right exact contravariant functor ${\rm
  Hom}_{H[\Lambda]}(\ \cdot\ , k^\vee_\bullet)$ as a composition of
  two right exact functors in three different ways as
  \begin{align*}
    {\rm Hom}_{H[\Lambda]}(\ \cdot\ ,  k^\vee_\bullet)
    \cong {\rm Hom}_k\left((\ \cdot\ )\eotimes{H}k\eotimes{k[\Lambda]}k_\bullet, k\right)
    \cong {\rm Hom}_k\left((\ \cdot\ )\eotimes{H[\Lambda]}H_\bullet\eotimes{H}k, k\right)
    \cong {\rm Hom}_H\left((\ \cdot\ )\eotimes{H[\Lambda]}H_\bullet, k\right)
  \end{align*}
  The result follows after using the fact that $k$ is a field,
  Proposition~\ref{CyclicHomology}, and  a Grothendieck spectral sequence
  argument.
\end{proof}

\begin{prop}\label{CoAlgebra}
  Let $Y_\bullet$ be a cocyclic $H$--module.  Then there are two
  spectral sequences with 
  \begin{align*}
   E^2_{p,q}=HC^p({\rm Ext}_H^q(k,Y_\bullet)) & & 
   'E^2_{p,q}={\rm Ext}_H^p(k,HC^q(Y_\bullet))  
  \end{align*}
  converging to ${\rm Ext}_{H[\Lambda]}^*(k_\bullet,Y_\bullet)$.
\end{prop}

\begin{proof}
  Note that we can write the left exact functor ${\rm
  Hom}_{H[\Lambda]}(k_\bullet, \ \cdot\ )$ as a composition of left
  exact functors in two different ways as
  \begin{align*}
    {\rm Hom}_{H[\Lambda]}(k_\bullet, \ \cdot\ )
    \cong {\rm Hom}_{k[\Lambda]}(k_\bullet, {\rm Hom}_H(k,\ \cdot\ ))
    \cong {\rm Hom}_H(k, {\rm Hom}_{H[\Lambda]}(H_\bullet,\ \cdot\ ))
  \end{align*}
  The result follows after using Proposition~\ref{CocyclicHomology}
  another Grothendieck spectral sequence argument, this time for
  composition of left exact functors.
\end{proof}

\begin{prop}\label{CoAlgebra2}
  There are two spectral sequences with
  \begin{align*}
    E_2^{p,q} = HC^p({\rm Tor}^H_q(Y^\vee_\bullet,k))  & & 
    'E_2^{p,q} = {\rm Ext}_H^p(HC_q(Y^\vee_\bullet),k)
  \end{align*}
  converging to ${\rm Ext}_{H[\Lambda]}^*(Y_\bullet,k_\bullet)$.
\end{prop}

\begin{proof}
  The proof of the first claim follows after observing the following
  sequence of isomorphism of functors
  \begin{align*}
    {\rm Hom}_{H[\Lambda]}(\ \cdot\ , k_\bullet)
    \cong {\rm Hom}_{H[\Lambda]}((\ \cdot\ )^\vee, k^\vee_\bullet)
    \cong {\rm Hom}_k\left((\ \cdot\ )^\vee\eotimes{H[\Lambda]}k_\bullet, k\right)
    \cong {\rm Hom}_k\left((\ \cdot\ )^\vee\eotimes{H}k\eotimes{k[\Lambda]}k_\bullet, k\right)
  \end{align*}
  by using a Grothendieck spectral sequence argument.  For the second
  claim we observe another sequence of isomorphisms
  \begin{align*}
    {\rm Hom}_{H[\Lambda]}(\ \cdot\ , k_\bullet)
    \cong {\rm Hom}_k\left((\ \cdot\ )^\vee\eotimes{H[\Lambda]}k_\bullet, k\right)
    \cong {\rm Hom}_H\left((\ \cdot\ )^\vee\eotimes{H[\Lambda]}H_\bullet, k\right)
  \end{align*}
  and we use another Grothendieck spectral sequence.
\end{proof}

Let $\C{A}$ be one of $\coalg{H}$, $H[\Lambda_\B{Z}]$, $H[\Lambda]$,
$H[\C{D}_\B{N}]$ or $H[\C{D}_C]$.

\begin{prop}
  There is a morphism of algebras of the form
  $\C{A}\xra{\lambda}H\otimes\C{A}$.  This yields two functors
  $L_\lambda:\lmod{H\otimes\C{A}}\xra{} \lmod{\C{A}}$ and $R_\lambda
  :\rmod{H\otimes\C{A}}\xra{} \rmod{\C{A}}$.
\end{prop}

\begin{proof}
  We will give the proof for $\coalg{H}$.  The proofs for the other
  cases are similar.  We define $\lambda$ on the generators by
  \begin{align*}
    \lambda(\partial^n_j) = 1_H\otimes\partial^n_j & &
    \lambda(\sigma^n_i) = 1_H\otimes\sigma^n_i & &
    \lambda(\tau_n^\ell) = 1_H\otimes\tau_n^\ell & & 
    \lambda(h) = h_{(1)}\otimes h_{(2)}
  \end{align*}
  for any $n,m\geq 0$, $\ell\in\B{Z}$ and $h\in H$.  Since $H$ is a
  bialgebra, we know that $H\xra{\Delta}H\otimes H$ is a morphism of
  $k$--algebras.  This proves the first assertion.  As for the second
  assertion, given a left $H\otimes\C{A}$--module $X_\bullet$ we
  define $L_\lambda(X_\bullet)=X_\bullet$ on the object level.  We define
  a left $\C{A}$--module structure on $L_\lambda(X_\bullet)$ by letting
  \begin{align*}
    \Psi(x):= (\Psi_{(-1)}\otimes\Psi_{(0)})(x)
  \end{align*}
  for any $x\in X_\bullet$ and for any $\Psi\in\C{A}$ where we use
  Sweedler's notation $\lambda(\Psi)=(\Psi_{(-1)}\otimes\Psi_{(0)})$.
\end{proof}

\begin{thm}\label{CohomologyOperations}
  For any left (or right) $\C{A}$--modules $X_\bullet$ and $Y_\bullet$
  the bivariant cohomology ${\rm Ext}_\C{A}^*(X_\bullet,Y_\bullet)$ is
  a graded ${\rm Ext}_H^* (k,k)$--module.
\end{thm}

\begin{proof}
  We have an action
  \begin{align*}
    {\rm Ext}_H^p(k,k)\otimes    
      {\rm Ext}_\C{A}^q(X_\bullet,Y_\bullet)\xra{}
      {\rm Ext}_{H\otimes\C{A}}^{p+q}(k\otimes X_\bullet,k\otimes Y_\bullet)
      \xra{L_\lambda^*}{\rm Ext}_\C{A}^{p+q}(X_\bullet,Y_\bullet)
  \end{align*}
  The action is associative since $\C{A}\xra{\lambda}H\otimes\C{A}$ is
  a coassociative coaction and $L_\lambda(k_\bullet\otimes X_\bullet)
  \cong X_\bullet$ for any $\C{A}$--module $X_\bullet$.
\end{proof}

\begin{prop}\label{BigBialgebra}
  $\C{A}$ is a non-unital, non-counital bialgebra.
\end{prop}

\begin{proof}
  We will give the proof for $\C{A}=\coalg{H}$.  The proofs for the
  other cases are similar.  We define a comultiplication $\Delta$ by
  defining it only on the generators.  We let $\Delta(h)=h_{(1)}\otimes
  h_{(2)}$ for any $h\in H$.  The rest of the generators are group-like, i.e.
  \begin{align*}
    \Delta(\partial^n_i)=\partial^n_i\otimes \partial^n_i && 
    \Delta(\sigma^n_j)=\sigma^n_j\otimes \sigma^n_j && 
    \Delta(\tau_n^\ell)= \tau_n^\ell\otimes \tau_n^\ell
  \end{align*}
  for all possible $m,n,i,j,\ell$.  It is routine to check that the
  comultiplication gives us a morphism of algebras
  $\C{A}\xra{\Delta}\C{A}\otimes\C{A}$.  
\end{proof}

\begin{thm}
  ${\rm Ext}_\C{A}^*(k_\bullet,k_\bullet)$ is a graded $k$--algebra
  and the bivariant cohomology groups ${\rm
  Ext}_\C{A}^*(X_\bullet,Y_\bullet)$ form a graded $k$--module over
  ${\rm Ext}_\C{A}^*(k_\bullet,k_\bullet)$.
\end{thm}

\begin{proof}
  The proof is very similar to the proof of
  Theorem~\ref{CohomologyOperations} after observing that $\C{A}$ is a
  bialgebra is equivalent to the fact that
  $\C{A}\xra{\Delta}\C{A}\otimes\C{A}$ is a coassociative coaction and that
  $L_\Delta(k_\bullet\otimes X_\bullet)$ is naturally isomorphic to
  $X_\bullet$ as $\C{A}$--modules.
\end{proof}

\begin{cor}\label{BottAlgebra}
  There is a spectral sequence of graded algebras with
  \begin{align*}
    E_2^{p,q} = {\rm  Ext}_H^p(k,k)\otimes HC^q(k_\bullet)  
  \end{align*}
  converging to ${\rm Ext}_{H[\Lambda]}^*(k_\bullet,k_\bullet)$.
\end{cor}

\begin{proof}
  We will use the spectral sequence in Proposition~\ref{CoAlgebra}.
  Note that ${\rm Ext}_{H[\Lambda]}^*(H_\bullet,k_\bullet)$ is the
  ordinary cyclic cohomology of $k$ viewed as a coalgebra.  Therefore
  it is the polynomial algebra $k[u]$ generated by an element of
  degree $2$ with the trivial $H$--action.
\end{proof}

  The bivariant theory we defined above is essentially an $H$--linear
  version of Connes' bivariant cyclic theory \cite{Connes:ExtFunctors, Connes:NonCommutativeGeometry}
  as studied by \cite{Nistor:BivariantChernConnes}.  One can also define a
  Jones-Kassel \cite{JonesKassel:BivariantCyclicTheory} variant of the
  cohomology we developed above as follows.  Instead of using the
  derived category of $H[\Lambda]$--modules, one can use the derived
  category of unbounded $H$--linear chain complexes.  For the
  positive cyclic theory we first chose an arbitrary projective
  resolution $\C{R}_{\bullet,*}$ of $k_\bullet$ (cocyclic $H$--module
  of $k$ considered as a $H$--module coalgebra) in the category of
  cocyclic $H$--modules and consider the differential graded
  $k$--module $X_\bullet\eotimes{H[\Lambda]}\C{R}_{\bullet,*}$ for a
  cyclic $H$--module $X_\bullet$ and ${\rm
  Hom}_{H[\Lambda]}(\C{R}_{\bullet,*},Y_\bullet)$ for a cocyclic
  $H$--module $Y_\bullet$.  These associations are functorial.  Denote
  these functors by $CP^{JK}_*(X_\bullet)$ and $CP_{JK}^*(Y_\bullet)$
  respectively.  As long as our (co)cyclic objects are H-(co)unital,
  one can replace $C^{JK}_*(X_\bullet)$ and $C_{JK}^*(Y_\bullet)$ by
  the corresponding $(b,B)$--complexes of $X_\bullet$ and $Y_\bullet$.  These
  complexes have an intrinsic degree $\pm 2$ endomorphism due to the
  fact that cohomology of cyclic groups have such an endomorphism.  This is Connes' 
  periodicity operator $S$.  Therefore,
  $C^{JK}_*(X_\bullet)$ and $C_{JK}^*(Y_\bullet)$ are differential
  graded $H[S]$--modules where $\deg(S)=\pm 2$ depending on whether we
  have a cyclic or cocyclic $H$--module.  Then for two right
  $H[\Lambda]$--modules $X_\bullet$ and $X_\bullet'$ and two left
  $H[\Lambda]$--modules $Y_\bullet$ and $Y_\bullet'$ we define two
  bivariant functors
  \begin{align*}
    HC^*_{JK}(X_\bullet,X_\bullet'):= &
      {\rm\bf Ext}_{H[S]}^*(CP^{JK}_*(X_\bullet), CP^{JK}_*(X_\bullet'))\\
    HC^*_{JK}(Y_\bullet,Y_\bullet'):= &
      {\rm\bf Ext}_{H[S]}^*(CP_{JK}^*(Y_\bullet), CP_{JK}^*(Y_\bullet'))
  \end{align*}
  where ${\rm\bf Ext}$ functors are the morphisms in the derived
  category of unbounded chain complexes of $H[S]$--modules.  Since the
  associations $X_\bullet\mapsto CP_*^{JK}(X_\bullet)$ and
  $Y_\bullet\mapsto CP^*_{JK}(Y_\bullet)$ are functorial, we have
  well-defined morphisms, called comparison morphisms, as
  \begin{align*}
    {\rm Ext}_{H[\Lambda]}^*(X_\bullet,X_\bullet') & \xra{CP_{JK}}
      {\rm\bf Ext}_{H[S]}^*(CP^{JK}_*(X_\bullet), CP^{JK}_*(X_\bullet'))\\
    {\rm Ext}_{H[\Lambda]}^*(Y_\bullet,Y_\bullet') & \xra{CP^{JK}}
      {\rm\bf Ext}_{H[S]}^*(CP_{JK}^*(Y_\bullet), CP_{JK}^*(Y_\bullet'))
  \end{align*}
  There are also negative and periodic versions of this bivariant
  cohomology theory which admit comparison morphisms, similar to the
  comparison morphisms we gave above, from our bivariant theory.  For
  the negative and periodic cyclic theories we replace
  $\C{R}_{\bullet,*}$ with two specific unbounded differential graded
  cocyclic $H[S]$--modules.

\section{Bivariant Hopf and bivariant equivariant cyclic cohomology of module algebras}
\label{ModuleAlgebra}

Recall that an algebra $A$ is called a right $H$--module algebra if
(i) $A$ is a right $H$--module and (ii) the multiplication map
$A\otimes A\to A$ is an $H$--module morphism, i.e. $(a_1 a_2)h =
(a_1h_{(1)}) (a_2h_{(2)})$ for any $h\in H$ and $a_1, a_2\in A$.  If
$A$ is unital, we also require $(1_A)h=\varepsilon(h)1_A$ for any
$h\in H$.

\begin{defn}
  Let $A$ be an $H$--module algebra and $M$ be a right-right
  $H$--module/comodule.  Define a right $\coalg{H}$--module
  $\B{T}_\bullet(A,M)$ by letting $\B{T}_n(A,M):= A^{\otimes
  n+1}\otimes M$ and
  \begin{align*}
    (a_0\otimes\cdots\otimes a_n\otimes m)\tau_n
    = & a_n m_{(1)}\otimes a_0\otimes\cdots\otimes a_{n-1}\otimes m_{(0)}\\
    (a_0\otimes\cdots\otimes a_n\otimes m)\partial^{n-1}_j
    = & \begin{cases}
        a_0 \otimes \cdots\otimes a_j a_{j+1}\otimes\cdots\otimes m 
                  & \text{ if } 0\leq j\leq n-1\\
        (a_n m_{(1)})a_0\otimes a_1\otimes\cdots\otimes a_{n-1}\otimes m_{(0)}
                  & \text{ if } j=n
	\end{cases}\\
    (a_0\otimes\cdots\otimes a_n\otimes m)\sigma^n_i
    = & a_0\otimes \cdots\otimes a_i\otimes 1_A\otimes a_{i+1}\otimes\cdots\otimes m \\
    (a_0\otimes\cdots\otimes a_n\otimes m)h
    = & (a_0 h_{(2)})\otimes \cdots\otimes (a_n h_{(n+2)})\otimes m h_{(1)}
  \end{align*}
   for $h\in H$, $0\leq j\leq n$ and $0\leq i\leq n,$ and 
  $a_0\otimes\cdots\otimes a_n\otimes m $ in $\B{T}_n(A,M)$.
\end{defn}

\begin{defn}
  Let $A$ and $M$ be as before.  We define $Q_\bullet(A,M)$ as the
  largest quotient of $\B{T}_\bullet(A,M)$ which is a cyclic
  $H$--module. Equivalently, we have
  \begin{align*}
    Q_\bullet(A,M) := \B{T}_\bullet(A,M)\eotimes{\coalg{H}}H[\Lambda]
  \end{align*}
\end{defn}

\begin{defn}
  Let $A$ and $A'$ be $H$--module algebras and $M$ and $M'$ be
  $H$--module/comodules.  We define the bivariant equivariant cyclic
  cohomology of the pair $(A,A')$ with coefficients in the pair
  $(M,M')$ as
  \begin{align*}
    HC_H^*(A,A';M,M'):={\rm Ext}_{H[\Lambda]}^*(Q_\bullet(A,M), Q_\bullet(A',M'))
  \end{align*}
\end{defn}

In the special case where $M=M'=k$ and $H=k$, our definition reduces
to Connes' bivariant cyclic theory
\cite{Connes:NonCommutativeGeometry, Nistor:BivariantChernConnes}:
\begin{align*}
  HC_k^*(A,A';k,k')\cong HC^*(A,A')
\end{align*}
It is also immediate from our definition that we have a graded
associative product, defined by composition
\begin{align*}
  HC_H^p(A,A';M,M')\otimes HC_H^q(A',A'';M',M'') \xra{}
  HC_H^{p+q}(A,A'';M,M'')
\end{align*}
In particular each bivariant group $HC_H^*(A,A;M,M)$ is a graded
associative algebra.

In Theorem~\ref{CohomologyOperations} we show that bivariant cyclic
cohomology of cyclic $H$--modules are graded modules over the graded
algebra ${\rm Ext}_H^*(k,k)$.  Specializing to our present case, we
obtain
\begin{prop} 
  The bivariant equivariant cyclic cohomology groups
  $HC_H^*(A,A';M,M')$ are graded modules over the graded $k$--algebra
  ${\rm Ext}_H^*(k,k)$.
\end{prop}

\begin{cor}
  Let $G$ be a group and let $A$ and $A'$ be $G$--algebras.  Then the
  bivariant equivariant cyclic cohomology groups
  $HC^*_{k[G]}(A,A';M,M')$ are graded modules over the group
  cohomology algebra $H^*(G,k)$ for any pair of $G$--module/comodules $M$, $M'$.
\end{cor}

\begin{cor}
  Let $\G{g}$ be a Lie algebra and let $A, $ and $A'$ be
  $U(\G{g})$--module algebras, i.e. $\G{g}$ acts on $A$ and $A'$ by
  derivations.  Then the bivariant equivariant cyclic cohomology
  groups $HC^*_{U(\G{g})}(A,A';M,M')$ are graded modules over the Lie
  algebra cohomology $H^*(\G{g},k)$ for any pair of
  $U(\G{g})$--module/comodules $M$, $M'$.
\end{cor}

Now we shift our attention to a different kind of bivariant cyclic
(co)homology of $H$--module algebras. 
\begin{defn}
  Given an $H$--module algebra $A$ and an $H$--module/comodule $M$ we
  define a cyclic $k$--module 
  \begin{align*}
    \C{C}_\bullet(A,M):= Q_\bullet(A,M)\eotimes{H}k
  \end{align*}
  With this definition at hand, we define Hopf cyclic homology and
  cohomology of the triple $(A,H,M)$ by
  \begin{align*}
    HC^\text{Hopf}_*(A,M):= {\rm Tor}^{k[\Lambda]}_*(\C{C}_\bullet(A,M),k_\bullet) & & 
    HC_\text{Hopf}^*(A,M):= {\rm Ext}_{k[\Lambda]}^*(\C{C}_\bullet(A,M),k^\vee_\bullet)
  \end{align*}
  Note that, if $H=k$ and $M=k$, this definition reduces to the ordinary
  cyclic (co)homology of algebras.  Note also that we can define another
  bivariant cyclic theory, called bivariant Hopf cyclic cohomology, of
  $H$--module algebras as
  \begin{align*}
    HC_{\rm Hopf}^*(A,A';M,M'):= {\rm
      Ext}_{k[\Lambda]}(\C{C}_\bullet(A,M),\C{C}_\bullet(A',C'))
  \end{align*}
  for any pair of $H$--module algebras $A$ and $A'$ and any pair of
  stable $H$--module/comodules $M$ and $M'$.
\end{defn}

\section{Bivariant Hopf and bivariant equivariant cyclic cohomology of module coalgebras}
\label{ModuleCoalgebra}

Recall that a coalgebra $C$ is called a left $H$--module coalgebra if
(i) $C$ is a left $H$--module and (ii) the comultiplication map
$C\xra{\Delta} C\otimes C$ is an $H$--module morphism, i.e. one has
$\Delta(h(c))= h_{(1)}(c_{(1)})\otimes h_{(2)}(c_{(2)})$ for any $h\in
H$ and $c\in C$.  Also, if $C$ is counital then we require
$\varepsilon(h(c))=\varepsilon(h)\varepsilon(c)$ for any $h\in H$ and
$c\in C$.

\begin{defn}
  Let $C$ be a left $H$--module coalgebra and $M$ be a left-left
  $H$--module/comodule.  We define a left $\coalg{H}$--module
  $\B{T}_\bullet(C,M)$ by  letting
  $$\B{T}_n(C,M)=C^{\otimes
  n+1}\otimes M$$ with the $\coalg{H}$--action defined as
  \begin{align*}
    \tau_n(c^0\otimes\cdots\otimes c^n\otimes m)
    = & c^1\otimes\cdots\otimes c^n\otimes m_{(-1)}c^0\otimes m_{(0)}\\
    \partial^n_j(c^0\otimes\cdots\otimes c^n\otimes m)
    = & \begin{cases}
        c^0\otimes \cdots\otimes c^j_{(1)}\otimes c^j_{(2)}\otimes\cdots\otimes m
             & \text{ if } 0\leq j\leq n\\
        c^0_{(2)}\otimes c^1\otimes\cdots\otimes c^n\otimes m_{(-1)}c^0_{(1)}\otimes m_{(0)})
             & \text{ if } j=n+1
	\end{cases}\\
    \sigma^{n-1}_i(c^0\otimes\cdots\otimes c^n\otimes m)
    = & \varepsilon(c^i)(c^0\otimes \cdots\otimes c^{i-1}\otimes c^{i+1}\otimes\cdots\otimes c^n\otimes m\\
    h(c^0\otimes\cdots\otimes c^n\otimes m)
    = & h_{(1)}(c_0)\otimes\cdots\otimes h_{(n+1)}(c_n)\otimes h_{(n+2)}m
  \end{align*}
  defined for $0\leq j\leq n+1$, $0\leq i\leq n-1$, $h\in H$ and
  $c^0\otimes\cdots\otimes c^n\otimes m$ from $\B{T}_n(C,M)$.
\end{defn}

\begin{defn}\label{HopfCyclicModule}
  Assume $C$, $H$,  and $M$ are as above.  We define two cocyclic $k$--modules
  \begin{align*}
    Q_\bullet(C,M):= H[\Lambda]\eotimes{\coalg{H}}\B{T}_\bullet(C,M) & & 
    \C{C}_\bullet(C,M):= k[\Lambda]\eotimes{\coalg{H}}\B{T}_\bullet(C,M)
  \end{align*}
  The latter can also be defined as $\C{C}_\bullet(C,M) =
  k\eotimes{H}Q_\bullet(C,M)$.
\end{defn}

\begin{defn}
  We define the bivariant Hopf cyclic cohomology groups of a pair of
  $H$--module coalgebras $C, C'$ and a pair of $H$--module/comodules
  $M, M'$ as
  \begin{align*}
    HC^*_{\rm Hopf}(C,C';M,M')
    = {\rm Ext}_{k[\Lambda]}^*(\C{C}_\bullet(C,M),\C{C}_\bullet(C',M'))
  \end{align*}
\end{defn}

Recall from \cite{Khalkhali:SaYDModules} that an arbitrary right-right
(resp. left-left) $H$--module/comodule $M$ is called stable if one has
$m_{(0)}m_{(1)}=m$ (resp. $m_{(-1)}m_{(0)}=m$) for any $m\in M$. When
$H$ is a Hopf algebra with an invertible antipode, $M$ is called an
anti-Yetter-Drinfeld module if
\begin{align*}
  (hm)_{(-1)}\otimes (hm)_{(0)} = h_{(1)}m_{(-1)}S^{-1}(h_{(3)})\otimes h_{(2)}m_{(0)}
\end{align*}
for any $h\in H$ and $m\in M$.  Following
\cite{Khalkhali:SaYDModules}, we will refer stable
anti-Yetter-Drinfeld modules as SAYD-modules.

When $H$ is a bialgebra, in \cite{Kaygun:BialgebraCyclicK} for an
$H$--module coalgebra $C$ and a stable $H$--module/comodule $M$, we
defined a para-cocyclic $H$--module $\B{PCM}_*(C,H,M)$ as the largest
quotient of $\B{T}_\bullet(C,M)$ which is a para-cocyclic $H$--module.
In other words $\B{PCM}_*(C,H,M):=
H[\Lambda_\B{Z}]\eotimes{\coalg{H}}\B{T}_\bullet(C,M)$.  Then we used
the product $k\eotimes{H}\B{PCM}_*(C,H,M)$ to compute the bialgebra
cyclic homology of the triple $(C,H,M)$.  It is also shown in
\cite{Kaygun:BialgebraCyclicK} that when $H$ is a Hopf algebra with an
invertible antipode and $M$ is an SAYD module this complex reduces to
the Hopf cyclic complex triple $(C,H,M)$
\cite{Khalkhali:HopfCyclicHomology}.

\begin{lem}\label{Special1}
  Let $C$ be an $H$--module coalgebra and $M$ be a stable
  $H$--module/comodule.  Then one has an isomorphism of cocyclic
  $k$--modules $k\eotimes{H}\B{PCM}_*(C,H,M) \cong
  \C{C}_\bullet(C,M)$.  In other words, $\C{C}_\bullet(C,M)$ is
  isomorphic to the cocyclic $k$--module which yields bialgebra cyclic
  homology in \cite{Kaygun:BialgebraCyclicK}.
\end{lem}

\begin{proof}
  Notice that $\tau_n^{n+1}(c^0\otimes\cdots\otimes c^n\otimes m) =
  m_{(-1)}(c^0\otimes\cdots\otimes c^n\otimes m_{(0)})$ for any $n\geq
  0$ and $(c^0\otimes\cdots\otimes c^n\otimes m)$ in $Q_\bullet(C,M)$.
\end{proof}

Now, Lemma~\ref{Special1} implies that if $C=H$ is a Hopf algebra and
$M$ is a SAYD module then $\C{C}_\bullet(H,M)$ isomorphic to the
cocyclic $k$--module of \cite{Khalkhali:HopfCyclicHomology} which is
used to compute Hopf cyclic cohomology.  Which in turn means if $H$ is
a Hopf algebra which admits a modular pair $(\sigma.\delta)$ in
involution \cite{ConnesMoscovici:HopfCyclicCohomology} and
$M=k_{(\sigma,\delta)}$ is the 1-dimensional $H$--module/comodule
associated with the modular pair $(\sigma,\delta)$ in involution, then
$\C{C}_\bullet(H,k_{(\sigma,\delta)})$ is the cocyclic module of
\cite{ConnesMoscovici:HopfCyclicCohomology,
ConnesMoscovici:HopfCyclicCohomologyII} which computes Hopf cyclic
cohomology of such a Hopf algebra.

\begin{thm}[{\bf Specialization theorem}]\label{Specialization}
  Let $C$ be an $H$--module coalgebra and $M$ be a stable
  $H$--module/comodule.  Bivariant Hopf cyclic cohomology groups
  $HC^*_{\rm Hopf}(k,C;k,M)$ are the Hopf cyclic cohomology groups of
  the triple $(C,H,M)$.  Bivariant Hopf cyclic cohomology groups
  $HC^*_{\rm Hopf}(C,k;M,k)$ are the cyclic cohomology groups of the
  dual cyclic module $\C{C}_\bullet(C,M)^\vee$.
\end{thm}

\begin{proof}
  Proof easily follows from Proposition~\ref{CocyclicHomology} and
  Lemma~\ref{Special1}.
\end{proof}

As in the case of ordinary bivariant cyclic cohomology, the bivariant
Hopf cyclic cohomology carry an additional structure.  The proof of
the following Proposition is the same as its ordinary counterpart in
cyclic homology of algebras \cite{Nistor:BivariantChernConnes,
Connes:ExtFunctors}.
\begin{prop}
  Bivariant Hopf cyclic cohomology groups $HC^*_{\rm Hopf}(C,C';M,M')$
  of a pair of $H$--module coalgebras $C, C'$ and a pair of stable
  $H$--module/comodules $M, M'$ are graded modules over the polynomial
  algebra $k[u] = {\rm Ext}_{k[\Lambda]}(k_\bullet,k_\bullet)$ where
  $\deg(u)=2$.
\end{prop}

Now we shift our attention to a different kind of a bivariant cyclic
(co)homology of $H$--module coalgebras.
\begin{defn}
  For a pair of $H$--module coalgebras $C$ and $C'$ and a pair of
  $H$--module/comodules $M$ and $M'$ we define bivariant equivariant 
  cyclic homology
  \begin{align*}
    HC_H^*(C,C';M,M'):={\rm Ext}_{H[\Lambda]}^*(Q_\bullet(C,M), Q_\bullet(C',M'))
  \end{align*}
\end{defn}

\begin{prop}\label{Spectral1}
  Let $C$ be an $H$--module coalgebra and $M$ be a
  $H$--module/comodule.  Then there is a spectral sequence with
  \begin{align*}
    E_2^{p,q}=HC^p({\rm Ext}_H^q(k,Q_\bullet(C,M)))
  \end{align*}
  converging to the bivariant equivariant cyclic cohomology groups
  $HC_H^*(k,C;k,M)$.  In particular, if $H$ is semi-simple then
  $HC_H^*(k,C;k,M)$ is the Hopf cyclic homology of the pair $(C,M)$.
\end{prop}

\begin{proof}
  Proposition~\ref{CoAlgebra}.
\end{proof}

\begin{prop}\label{Spectral2}
  Let $C$ be an $H$--module coalgebra and $M$ be a
  $H$--module/comodule.  Then there is a spectral sequence with
  \begin{align*}
    E^2_{p,q}=HC^p({\rm Tor}^H_q(Q_\bullet(C,M)^\vee,k))
  \end{align*}
  converging to the bivariant equivariant cyclic cohomology groups
  $HC_H^*(C,k;M,k)$.  In particular for $q=0$ we see that $E^2_{p,0}$
  is the dual Hopf cyclic cohomology of the triple $(C,H,M)$.
\end{prop}

\begin{proof}
  We use Proposition~\ref{CoAlgebra2} to get the spectral sequence.
  Second assertion follows from Lemma~\ref{Special1}.
\end{proof}

The following Proposition is an immediate consequence of the
Theorem~\ref{CohomologyOperations}.
\begin{prop}
  The bivariant equivariant cohomology groups $HC^*_H(C,C';M,M')$ are
  graded modules over the graded algebra ${\rm Ext}_H^*(k,k)$.
\end{prop}

\begin{cor}
  Let $G$ be a group and let $C$ be a $G$--coalgebra.  Then the
  bivariant equivariant cyclic cohomology modules
  $HC^*_{k[G]}(C,C';M,M')$ are graded modules over the group
  cohomology $H^*(G,k)$ for any pair $C,C'$ of $k[G]$--module
  coalgebras and any pair $M,M'$ of $k[G]$--module/comodules.
\end{cor}

\begin{cor}
  Let $\G{g}$ be a Lie algebra and let $C$ be a $U(\G{g})$--module
  coalgebra.  Then the bivariant equivariant cyclic cohomology modules
  $HC^*_{U(\G{g})}(C,C';M,M')$ are graded modules over the Lie algebra
  cohomology $H^*(\G{g},k)$ for any pair $C,C'$ of $U(\G{g})$--module
  coalgebras and any pair $M,M'$ of $U(\G{g})$--module/comodules.
\end{cor}
 
\section{Cyclic (co)homology of $H$--categories and Morita invariance}
\label{MoritaInvariance}

Let  $H$ be  an arbitrary Hopf algebra.  Our goal in this section is
to extend the formalism of Hopf cyclic cohomology to a certain class
of additive categories that we call $H$--categories. We show that if
an $H$--category $\C{C}$ is separated over a subcategory $\C{E}$, the
Hopf cyclic and equivariant cyclic (co)homology groups of $\C{C}$ and
$\C{E}$ are isomorphic.  From this one easily derives Morita
invariance theorems for Hopf cyclic and equivariant cyclic
(co)homology theories we developed previous sections.

\begin{defn}
  A small $k$--linear category $\C{C}$ is called an $H$--category if
  (i) for every pair of objects $X$ and $Y$ the set of morphisms ${\rm
  Hom}_\C{C}(X,Y)$ is a left $H$--module, (ii)
  $h(id_X)=\varepsilon(h)id_X$ for any $X\in Ob(\C{C})$ and $h\in H$,
  and (iii) the composition of morphisms is $H$--equivariant, i.e.
  \begin{align*}
    h(fg)=h_{(1)}(f) h_{(2)}(g)
  \end{align*}
  for every pair $Z\xla{f}Y$ and $Y\xla{g}X$ of composable morphisms.  A
  functor of $H$--categories $F:\C{C}\xra{}\C{C}'$ is called an
  $H$--equivariant functor if the structure morphisms
  \begin{align*}
    F_{X,Y}:{\rm Hom}_\C{C}(X,Y)\xra{}{\rm Hom}_{\C{C}'}(F(X),F(Y))
  \end{align*}
  are morphisms of $H$--modules, i.e. $F_{X,Y}(h(f))=h(F_{X,Y}(f))$
  for any pair of objects $X$ and $Y$ in $\C{C}$ and for any $h\in H$
  and $f\in{\rm Hom}_\C{C}(X,Y)$.
\end{defn}

\begin{defn}
  A morphism $X\xla{f}Y$ in an $H$--category $\C{C}$ is called
  $H$--invariant if $h(f)=\varepsilon(h)f$ for any $h\in H$.  The
  subcategory of $H$-invariant morphisms of $\C{C}$ is denoted by
  ${}^H\C{C}$.
\end{defn}

Our main source of examples of $H$--categories will come from
$H$--equivariant modules over $H$--module algebras.
\begin{defn}
  Let $A$ be a unital $H$--module algebra.  A left $A$--module $P$ is
  called an $H$--equivariant $A$--module if (i) $P$ is an $H$--module
  and (ii) the module action $A\otimes P \to P$ is a morphism of
  $H$--modules i.e.
  \begin{align*}
    h(a x) = h_{(1)}(a) h_{(2)}(x)
  \end{align*}
  for any $a \in A$, $h\in H$ and $x\in P$. Similarly, one can define
   $H$--equivariant left $A$--modules and $H$--equivariant
  $A$--bimodules.
\end{defn}

\begin{exm}
  Let $H$ be a Hopf algebra and $A$ be an arbitrary $H$--module
  algebra.  Let $\emod{H}{A}$ be the category of all finitely
  generated $H$--equivariant right $A$--modules with all morphisms of
  right $A$--modules between them.  For a pair of objects $X$ and $Y$
  we define an $H$--action on ${\rm Hom}_A(X,Y)$ by letting
  $h(f)(x)=h_{(1)}f(S(h_{(2)})x)$ for any $x\in X$ and $f\in {\rm
  Hom}_k(X,Y)$.  It is easy to check that $h(f)$ is still a morphism
  of right $A$--modules:
  \begin{align*}
    (hf)(xa) = & h_{(1)}f(S(h_{(2)})(xa)) = h_{(1)}f(S(h_{(3)})(x)S(h_{(2)})(a))\\
       = & h_{(1)}f(S(h_{(4)})(x)) h_{(2)}S(h_{(3)})(a) = (hf)(x) a
  \end{align*}
  for any $x\in X$, $a\in A$ and $h\in H$.  So $\emod{H}{A}$ is an
  $H$--category.  Note that ${}^H\emod{H}{A}$ is the subcategory of
  finitely generated $H$--equivariant $A$--modules with $H$--linear
  and $A$--linear morphisms.
\end{exm}

\begin{exm}
  Let $H$ and $A$ be as before.  One can also consider $\efree{H}{A}$
  the category of finitely generated free $A$--modules, which are
  $H$--equivariant $A$--modules automatically.  There is also
  $*_H^A$ the category which consists of only one object $A$
  considered as an $H$--equivariant right $A$--module.
\end{exm}

\begin{exm}
  A projective $A$--module $P$ is called $H$--equivariantly projective
  if (i) there exists another $A$--module $Q$ such that $F=P\oplus Q$
  is a free $A$--module (with the inherent $H$--equivariant
  $A$--module structure) and (ii) the canonical epimorphism
  $F\xra{\pi_P}P$ is $H$--invariant.  The subcategory of
  $\emod{H}{A}$ which consists of finitely generated
  $H$--equivariantly projective right $A$--modules is denoted by
  $\eproj{H}{A}$.
\end{exm}

\begin{defn}
  Let $\C{C}$ be an $H$--category and $M$ be a stable
  $H$--module/comodule.  We define a pseudo-para-cyclic module (i.e. a
  right $\coalg{H}$--module) $\B{T}_\bullet(\C{C},M)$ by letting
  \begin{align*}
    \B{T}_n(\C{C},M) 
    = & \bigoplus_{X_0,\ldots,X_n} {\rm Hom}_\C{C}(X_1,X_0)\otimes
           \cdots\otimes{\rm Hom}_\C{C}(X_n,X_{n-1})\otimes {\rm Hom}_\C{C}(X_0,X_n)
           \otimes M
  \end{align*}
  We define structure morphisms as follows:
  \begin{align*}
    (f_0\otimes\cdots\otimes f_n\otimes m)\tau_n
    = & m_{(-1)}f_n\otimes f_0\otimes\cdots\otimes f_{n-1}\otimes m_{(0)}\\
    (f_0\otimes\cdots\otimes f_n\otimes m)\partial^{n-1}_j
    = & \begin{cases}
        f_0\otimes \cdots\otimes f_j f_{j+1}\otimes\cdots\otimes f_n\otimes m 
             & \text{ if } 0\leq j\leq n-1\\
        (m_{(-1)}f_n)f_0\otimes f_1\otimes\cdots\otimes f_{n-1}\otimes m_{(0)}
             & \text{ if } j=n
	\end{cases}\\
    (f_0\otimes\cdots\otimes f_n\otimes m)\sigma^n_j
    = & \begin{cases}
        f_0\otimes\cdots\otimes f_j\otimes id_{X_{j+1}}\otimes f_{j+1}\otimes\cdots
         \otimes f_n\otimes m & \text{ if } 0\leq j\leq n-1\\
	f_0\otimes\cdots\otimes f_n\otimes id_{X_0}\otimes m
	                       & \text{ if } j=n
	\end{cases}\\
    (f_0\otimes\cdots\otimes f_n\otimes m)h
    = & S^{-1}(h_{(n+2)})(f_0)\otimes\cdots\otimes S^{-1}(h_{(2)})(f_n)\otimes S^{-1}(h_{(1)})(m)
  \end{align*}
  for any $n\geq 0$, $0\leq j\leq n$, $h \in H$,  and 
  $f_0\otimes\cdots\otimes f_n\otimes m \in\B{T}_n(\C{C},M)$.
\end{defn}

\begin{defn}
  For an $H$--category $\C{C}$ and an stable $H$--module/comodule $M$
  we define $Q_\bullet(\C{C},M)$ as the largest quotient of
  $\B{T}_\bullet(\C{C},M)$ which is a cyclic $H$--module.  In other
  words we divide $\B{T}_\bullet(\C{C},M)$ by the graded
  $H$--submodule and a cyclic $k$--module generated by elements of the
  form
  \begin{align*}
    (\Psi)[x,\tau_n^i]+(\Phi)(\tau_n^{n+1}-\tau_n^0)
  \end{align*}
  for $x\in H$, $n,i\in\B{N}$ and $\Psi,\Phi\in\B{T}_n(\C{C},M)$.  Equivalently we have 
  \begin{align*}
    Q_\bullet(\C{C},M):=\B{T}_\bullet(\C{C},M)\eotimes{\coalg{H}}H[\Lambda]
  \end{align*}
\end{defn}

Assume $\C{A}_*$ and $\C{B}_*$ are two pre-simplicial $H$--modules
(i.e. right $H[\C{D}]$--modules), which can also be considered as
differential graded $H$--modules.  Assume that $f,g:\C{A}_*\to\C{B}_*$
are morphisms of pre-simplicial $H$--modules. Recall that if there is
a set of $H$--module morphisms $h_j:\C{A}_n\to \C{B}_{n+1}$ defined
for $0\leq j\leq n$ satisfying
\begin{align*}
  h_j\partial_i = & \partial_i h_{j+1},\ \  \text{ if } i\leq j &
  h_j \partial_i = & \partial_{i+1} h_j,\ \ \text{ if } i\geq j+1 &  
  \partial_i h_i = & \partial_i h_{i-1},\ \ \text{ if } i\geq 1
\end{align*}
where $f_n=\partial_0 h_0$ and $g_n=\partial_{n+1} h_n$ for any $n\geq
0$, then $s_n = \sum_{i=0}^n (-1)^i h_i$ defines a chain homotopy
between the morphism of differential graded $H$--modules $f_*$ and
$g_*$.

For any (co)cyclic $H$--module $X_\bullet$ (i.e. a right and
resp. left $H[\Lambda]$--module), we will use $X_\bullet^\C{D}$ to
denote the same (co)cyclic $H$--module viewed as a pre-(co)simplicial
module (i.e. right and resp. left $H[\C{D}]$--module).  In other
words, $X_\bullet^\C{D}$ stands for the restriction ${\rm
Res}^{H[\Lambda]}_{H[\C{D}]}(X_\bullet)$.

\begin{lem}\label{PreSimplicialHomotopy}
  Fix a stable $H$--module/comodule $M$.  If two $H$--categories
  $\C{C}$ and $\C{C}'$ are $H$--equivariantly equivalent then
  $Q_\bullet^\C{D}(\C{C},M)$ and $Q_\bullet^\C{D}(\C{C}',M)$ are
  $H$--equivariantly homotopy equivalent.
\end{lem}

\begin{proof}
  Assume we have a pair of $H$--equivariant functors
  $F:\C{C}\xra{}\C{C}'$ and $G:\C{C}'\xra{}\C{C}$ such that
  $id_\C{C}\simeq GF$ and $id_{\C{C}'} \simeq FG$ $H$--equivariantly,
  i.e. the natural transformations are $H$ invariant.  Call these
  natural transformations $h$ and $h'$ respectively.  We will show
  that $Q_\bullet^\C{D}(GF,M)\simeq id$ on
  $Q_\bullet^\C{D}(\C{C},M)$.  The proof for
  $Q_\bullet^\C{D}(FG,M)\simeq id$ on $Q_\bullet^\C{D}(\C{C}',M)$ is
  similar.  Define
  \begin{align*}
    h_j(f_0\otimes & \cdots\otimes f_n\otimes m)\\
    = & \begin{cases}
        f_0 h_{X_1}^{-1}\otimes GF(f_1)\otimes\cdots\otimes GF(x_j)\otimes 
        h_{X_{j+1}}\otimes f_{j+1}\otimes\cdots\otimes  f_n\otimes m
             & \text{ if } 0\leq j\leq n-1\\
        f_0 h_{X_1}^{-1}\otimes GF(f_1)\otimes\cdots\otimes GF(f_n)\otimes 
        h_{X_0}\otimes m
             & \text{ if } j=n
	\end{cases}
  \end{align*}
  for any $(f_0\otimes\cdots\otimes f_n\otimes m)$ in $Q_n(\C{C},M)$.
  Note that $\partial_0 h_0$ is the identity on
  $Q_\bullet^\C{D}(id_\C{C},M)$ and $\partial_{n+1} h_n =
  Q_\bullet^\C{D}(GF,M)$.  We leave the verification of the
  pre-simplicial homotopy identities to the reader.  This proves
  $Q_\bullet^\C{D}(GF,M)$ is homotopic to
  $Q_\bullet^\C{D}(id_\C{C},M)$.  The proof for $FG$ is similar.
\end{proof}

\begin{defn}
  Let $\C{C}$ be an $H$--category and let $\C{E}$ be a full
  $H$--subcategory of $\C{C}$.  The category $\C{C}$ is called
  separated over $\C{E}$ iff for every object $X$ in $\C{C}$ there
  exists a natural number $n\geq 1$ and a finite set of objects
  $\{X_1,\ldots,X_n\}$ in $\C{E}$ such that we have $H$--invariant
  morphisms $X_i\xla{u_i}X$ and $X\xla{v_j}X_j$ in $\C{C}$ which
  satisfy $u_iv_i=id_{X_i}$ and $\sum_{i=1}^n v_iu_i=id_X$.
\end{defn}

\begin{prop}\label{SeperatedChain}
  $\eproj{H}{A}$ is separated over $\efree{H}{A}$ which in turn is
  separated over $*_H^A$.
\end{prop}

\begin{lem}\label{Seperation}
  Let $\C{C}$ be an $H$--category separated over an $H$--subcategory
  $\C{C}$.  Then the natural inclusion of functors $\C{E}\xra{}\C{C}$
  induces a homotopy equivalence of pre-simplicial $H$--modules
  $Q_\bullet^\C{D}(\C{E},M)\xra{}Q_\bullet^\C{D}(\C{C},M)$ for any
  stable $H$--module/comodule $M$.
\end{lem}

\begin{proof}
  For every object $X$ in $\C{C}$ fix a finite set of objects
  $\{X_1,\ldots,X_{n(X)}\}$ and $H$--invariant morphisms
  $X_i\xla{u_i(X)}X$ and $X\xla{v_j(X)}X_j$ such that the choices for
  every object $X$ in $\C{E}$ are $\{X\}$ with $u_1(X)=id_X=v_1(X)$.
  Define a morphism of graded $H$--modules
  $Q_\bullet(\C{C},M)\xra{E_\bullet}Q_\bullet(\C{E},M)$ by letting
  \begin{align*}
    E_n(f_0\otimes & \cdots\otimes f_n\otimes m)\\
    = & \sum_{i_0,\ldots,i_n}
        u_{i_0}(X_0)f_0v_{i_1}(X_1)\otimes u_{i_1}(X_1)f_1v_{i_2}(X_2)\otimes
        \cdots \\ 
      & \hspace{1cm} \otimes u_{i_{n-1}}(X_{n-1})f_{n-1}v_{i_n(X_n)}\otimes
        u_{i_n}(X_n)f_nv_{i_0}(X_0)\otimes m
  \end{align*}
  for any $f_0\otimes\cdots\otimes f_n\otimes m$ in
  $Q_\bullet(\C{C},M)$.  One can easily check that $E_\bullet$ is a
  morphism of cyclic $H$--modules and the composition
  $Q_\bullet(\C{E},M)\xra{i_\bullet}Q_\bullet(\C{C},M)\xra{E_\bullet}
  Q_\bullet(\C{E},M)$ is the identity.  In order to prove that
  $i_\bullet E_\bullet$ and $id_\bullet$ are homotopic on
  $Q_\bullet^\C{D}(\C{C},M)$, we must furnish a pre-simplicial
  homotopy.  For any $n\geq 0$ and $(f_0\otimes \cdots\otimes
  f_n\otimes m)$ in $Q_n(\C{C},M)$ we define
  \begin{align*}
    h_j(f_0\otimes & \cdots\otimes f_n\otimes m)\\
    = & \sum_{i_0,i_{j+1},\ldots,i_n}
        u_{i_0}(X_0)f_0\otimes f_1\otimes\cdots\otimes f_j
        \otimes v_{i_{j+1}}(X_{j+1})\otimes 
        u_{i_{j+1}}(X_{j+1})f_{j+1}v_{i_{j+2}}(X_{j+2})\otimes\cdots\\
      & \hspace{2cm}
        \otimes u_{i_n}(X_n)f_n v_{i_0}(X_0) \otimes m
  \end{align*}
  for $j\geq 0$.  For $j=n$ we let
  \begin{align*}
    h_n(f_0\otimes & \cdots\otimes f_n\otimes m)
    = \sum_{i_0} u_{i_0}(X_0)f_0\otimes f_1\otimes\cdots\otimes f_n 
        \otimes v_{i_0}(X_0)\otimes m
  \end{align*}
  One can easily see that $\partial_0 h_0 = i_\bullet E_\bullet$ and
  $\partial_{n+1}h_n = id_\bullet$.  We leave the verification of the
  simplicial homotopy identities to the reader.
\end{proof}

\begin{defn}
  Two $H$--module algebras $A$ and $A'$ are said to be
  $H$--equivariantly Morita equivalent if there is an $H$--equivariant
  $A$--$A'$--bimodule $P$ and an $H$--equivariant $A'$--$A$--bimodule
  $Q$ such that $P\eotimes{A'}Q\cong A$ as $H$--equivariant
  $A$--$A$--bimodules and $Q\eotimes{A}P\cong A'$ as $H$--equivariant
  $A'$--$A'$--bimodules.
\end{defn}

\begin{thm}[{\bf Equivariant Morita invariance for module algebras}]\label{Morita}
  Assume $A$ and $A'$ are two unital $H$--module algebras which are
  $H$--equivariantly Morita equivalent.  Then Hopf cyclic (co)homology
  and equivariant cyclic (co)homology of the pairs $(A,M)$ and
  $(A',M)$ are isomorphic.
\end{thm}

\begin{proof}
  Lemma~\ref{PreSimplicialHomotopy} implies that pre-simplicial
  $H$--modules $Q_\bullet^\C{D}(\eproj{H}{A},M)$ and
  $Q_\bullet^\C{D}(\eproj{H}{A'},M)$ are homotopy equivalent.  The
  inclusion $*_H^A\xra{}\eproj{H}{A}$ induces a homotopy equivalence
  $Q_\bullet^\C{D}(*_H^A,M)\xra{} Q_\bullet^\C{D}(\eproj{H}{A},M)$ by
  Proposition~\ref{SeperatedChain} and Lemma~\ref{Seperation}.
  Observe also that $Q_\bullet(*_H^A,M)$ is isomorphic to
  $Q_\bullet(A,M)$ as $H[\Lambda]$--modules.  The result follows from
  the SBI-sequence in cyclic (co)homology and the 5-Lemma.
\end{proof}

\begin{exm}
  Let $C$ be a left $H$--module coalgebra.  A right $C$--comodule $X$
  is called an $H$--equivariant right $C$--comodule if the coaction
  $X\xra{\rho_X}X\otimes C$ is a left $H$--module morphism.
  Explicitly, for any $x\in X$ and $h\in H$ one has
  \begin{align*}
    (hx)_{(0)}\otimes (hx)_{(1)} = h_{(1)}x_{(0)}\otimes h_{(2)}x_{(1)}
  \end{align*}
  A morphism of $C$--comodules $X\xra{f}Y$ is said to be
  $H$--equivariant if $f$ is also $H$--linear.  An arbitrary (left)
  $H$--equivariant (right) $C$--comodule $F$ is said to be free if $F$
  is isomorphic to a direct sum $C^{\oplus I}$ with its canonical
  $H$--structure where $I$ is an arbitrary (not necessarily finite)
  index set.  If $I$ is finite, then $F$ is called finitely generated.
  An arbitrary $H$--equivariant $C$--comodule $P$ is said to be
  projective if there exists another $H$--equivariant $C$--comodule
  $Q$ such that (i) $F := P\oplus Q$ is a free $H$--equivariant
  $C$--comodule (ii) the canonical epimorphism $F\xra{\pi}C$ is
  $H$--invariant.  Let $\emod{H}{C}$ denote the category of all
  finitely generated $H$--equivariant $C$--comodules with not
  necessarily $H$--equivariant but all $C$--comodule morphisms.  We
  will use $\eproj{H}{C}$ and $\efree{H}{C}$ to denote the
  subcategories of finitely generated projective and finitely
  generated free $H$--equivariant $C$--comodules respectively.  There
  is also $*_H^C$ the subcategory of $\efree{H}{C}$ on one single
  object $C$ considered as an $H$--equivariant $C$--comodule via its
  comultiplication.

  The category $\emod{H}{C}$ and all of its subcategories we mentioned
  above carry an $H$--category structure defined as follows: given
  $X\xra{f}Y$ an arbitrary $C$--comodule morphism between two
  $H$--equivariant $C$--comodules and $h\in H$ we define
  \begin{align*}
    (hf)(x) = h_{(1)}f(S(h_{(2)})x)
  \end{align*}
  for any $x\in X$.  We must check that $hf$ is still a $C$--comodule
  morphism.  Recall that $f$ is a $C$--comodule morphism iff
  $f(x)_{(0)}\otimes f(x)_{(1)} = f(x_{(0)})\otimes x_{(1)}$ for any
  $x\in X$.  We also observe that since both $X$ and $Y$ are both
  $H$--equivariant $C$--comodules we have
  \begin{align*}
    ((hf)(x))_{(0)}\otimes ((hf)(x))_{(1)} 
    = & h_{(1)(1)}f(S(h_{(2)})x)_{(0)}\otimes  h_{(1)(2)}f(S(h_{(2)})x)_{(1)}\\
    = & h_{(1)}f(S(h_{(4)})x_{(0)})\otimes h_{(2)}S(h_{(3)})x_{(1)}\\
    = & (hf)(x_{(0)})\otimes x_{(1)}
  \end{align*}
  for any $x\in X$ and $h\in H$ as we wanted to show.

  As before, one can show that $\eproj{H}{C}$ is separated over
  $\efree{H}{C}$ which in turn is separated over $*_H^C$.  Thus Hopf
  and equivariant cyclic (co)homology of $Q_\bullet(\eproj{H}{C},M)$
  are the same as the Hopf and equivariant cyclic cohomology of the
  $H$--module algebra ${\rm Hom}_C(C,C)$ respectively for any stable
  $H$--module/comodule $M$.  Now let us identify the $H$--module
  algebra ${\rm Hom}_C(C,C)$.
\end{exm}

For a left $H$--module algebra $A$, the opposite $H$--module algebra
$A^{op}$ is a right $H$--module algebra with $x^{op} h:=
(S^{-1}(h)x)^{op}$ for any $x\in A$ and $h\in H$.  Note that
\begin{align*}
  (x^{op} y^{op}) h 
  = (S^{-1}(h)(yx))^{op} = (S^{-1}(h_{(2)})y S^{-1}(h_{(1)})x)^{op}
  = (x h_{(1)})^{op} (y h_{(2)})^{op}
\end{align*}
for any $x^{op},y^{op}\in A^{op}$ and $h\in H$.

\begin{prop}
  Let $C$ be a counital $H$--module coalgebra.  The $k$--linear dual
  $C^*$ is an $H$--module algebra with the convolution product.
  Moreover, $C^*$ is isomorphic to ${\rm Hom}_C(C,C)^{op}$ the
  opposite $H$--module algebra of $C$--comodule endomorphisms of $C$
  as $H$--module algebras.
\end{prop}

\begin{proof}
  Take two arbitrary elements $\delta,\mu\in C^*={\rm Hom}_k(C,k)$ and
  define
  \begin{align*}
    (\delta*\mu)(c) = \delta(c_{(1)})\mu(c_{(2)})
  \end{align*}
  for any $c\in C$.  One can easily check that $*$ is an associative
  product since $C$ is an coassociative coalgebra.  Note also that the
  counit $\varepsilon$ of $C$ is the unit of this convolution algebra.
  The right $H$--structure is given by $(\delta h)(c) = \delta(h c)$
  for any $c\in C$ and $h\in H$.  Now we can check that
  \begin{align*}
    ((\delta*\mu)h)(c) 
    = (\delta*\mu)(h c)
    = \delta(h_{(1)}c_{(1)})\mu( h_{(2)} c_{(2)})
    = ((\delta h_{(1)})*(\mu h_{(2)}))(c)
  \end{align*}
  for any $h\in H$ and $c\in C$.  This proves the first assertion.
  For the second assertion, observe that ${\rm Hom}_C(C,C)^{op}$ is a
  right $H$--module algebra since the action of $H$ on $*_H^C$ was on
  the left.  We define two $k$--linear morphisms ${\rm
  Hom}_C(C,C)^{op}\xra{\Phi}C^*$ and $C^*\xra{\Psi}{\rm
  Hom}_C(C,C)^{op}$ by
  \begin{align*}
    \Phi(u^{op})(c) := \varepsilon(u^{op}(c)) & & 
    \Psi(\delta)(c) := \delta(c_{(1)})c_{(2)}
  \end{align*}
  for any $u^{op}\in {\rm Hom}_C(C,C)^{op}$, $\delta\in C^*$ and $c\in
  C$.  It is easy to see that $\Phi(u^{op})$ is in $C^*$ for any
  $u^{op}\in {\rm Hom}_C(C,C)^{op}$.  On the other hand, for
  $\delta\in C^*$ we check
  \begin{align*}
    (\Psi(\delta)(c))_{(0)}\otimes (\Psi(\delta)(c))_{(1)}
    = \delta(c_{(1)})c_{(2)}\otimes c_{(3)} 
    = (\Psi(\delta)(c_{(1)}))\otimes c_{(2)}
  \end{align*}
  for any $c\in C$, i.e. $\Psi(\delta)$ lies in ${\rm
  Hom}_C(C,C)^{op}$.  Notice also that $\Psi$ and $\Phi$ are inverses
  of each other since
  \begin{align*}
    (\Psi(\Phi(u^{op})))(c)
     = & \varepsilon(u^{op}(c_{(1)}))c_{(2)}
     = \varepsilon(u^{op}(c)_{(1)})u^{op}(c)_{(2)} = u^{op}(c)\\
    (\Phi(\Psi(\delta)))(c)
     = & \varepsilon(\delta(c_{(1)})c_{(2)}) = \delta(c)
  \end{align*}
  for any $u^{op}\in {\rm Hom}_C(C,C)^{op}$, $\delta\in C^*$ and $c\in C$.
  For the compatibility of the $H$--structures we see that
  \begin{align*}
    \Phi(u^{op} h)(c) = \varepsilon(u^{op}(hc))
     = \varepsilon(S^{-1}(h_{(2)}) u^{op} (h_{(1)}c)) = \varepsilon(u^{op}(hc))
     = (\Phi(u^{op})h)(c)
  \end{align*}
  for any $u^{op}\in {\rm Hom}_C(C,C)^{op}$, $h\in H$ and $c\in C$.
  Finally we check
  \begin{align*}
    \Psi(\delta*\mu)(c) = (\delta*\mu)(c_{(1)})c_{(2)}
      = \delta(c_{(1)})\mu(c_{(2)})c_{(3)}
      = \Psi(\mu)(\delta(c_{(1)})c_{(2)}) = (\Psi(\mu)\circ\Psi(\delta))(c)
  \end{align*}
  where $\circ$ denotes composition of morphisms.  However, recall
  that in the opposite algebra the multiplication is the opposite
  composition.  Therefore
  \begin{align*}
    \Psi(\delta*\mu)=\Psi(\delta)^{op}\Psi(\mu)^{op}
  \end{align*}
  for any $\delta,\mu\in C^*$ as we wanted to show.
\end{proof}

\begin{thm}[{\bf Morita invariance for module coalgebras}]
  Let $M$ be an arbitrary stable $H$--module/ comodule.  Let $C$ and
  $C'$ be two $H$--module coalgebras such that the categories
  $\eproj{H}{C}$ and $\eproj{H}{C'}$ are $H$--equivariantly
  equivalent.  Then $H$--module algebras $C^*$ and $C'^*$ have
  isomorphic Hopf and equivariant cyclic (co)homologies.
\end{thm}

\begin{rem}
  The functor $(\ \cdot\ )^*:= {\rm Hom}_k(\ \cdot\ , k)$ sends an
  $H$--equivariant $C$--comodule to an $H$--equivariant $C^*$--module.
  However, the failure of $(\ \cdot\ )^*$ being an equivalence forbids
  us to translate the problem of equivariant Morita equivalence of two
  $H$--module coalgebras $C$ and $C'$ purely in terms of equivariant
  Morita equivalence of dual $H$--module algebras $C^*$ and $(C')^*$.
  Although it is true that if $C^*$ and $(C')^*$ are
  $H$--equivariantly Morita equivalent then $C^*$ and $(C')^*$ has
  isomorphic Hopf and equivariant cyclic cohomology, requiring $C$ and
  $C'$ to be equivariantly Morita equivalent is much weaker.
\end{rem}
  

\end{document}